\documentclass[12pt]{amsart}
\usepackage{latexsym}
\usepackage{amsmath,amsfonts,amssymb,amsthm}
\usepackage{graphicx}
\usepackage{color}
\usepackage{xypic}
\xyoption{all}
\DeclareMathOperator{\Pol}{Pol} %
\newcommand{\kk}{\K} %
\def\poldiv{\mathcal{D}} %
\DeclareMathOperator{\CaDiv}{CaDiv} %
\DeclareMathOperator{\tail}{tail} %
\DeclareMathOperator{\polf}{f} %
\newcommand{\Mb}{M^{\prime}} %
\newcommand{\Nb}{N^{\prime}} %
\newcommand{\Sb}{\Sigma^{\prime}}
\newcommand{\chQ}{/\hspace{-0.3em}/^{\mbox{\tiny ch}}} %
\newcommand{\toric}[1]{\mathbb T \mathbb V({#1})} %
\newcommand{\torusO}{T}
\newcommand{\torus}{\til{T}}
\newcommand{\kegel}{\sigma}
\DeclareMathOperator{\second}{Sec} %
\newcommand{\normal}{\mathcal N}
\newcommand{\kG}{\Gamma} %
\DeclareMathOperator{\orb}{orb} %
\DeclareMathOperator{\cone}{cone} %
\newcommand{\dcup}{\sqcup} %
\DeclareMathOperator{\face}{face} %
\DeclareMathOperator{\loc}{loc} %
\newcommand{\Ngrass}{\til{N}} %
\newcommand{\fanmap}{\pi'}
\newcommand{\umunu}{u^{\mu\nu}}
\newcommand{\Dmunu}{E^{\mu\nu}}
\newcommand{\MOb}{M''}
\newcommand{\MOm}{\ker\ku{1}}
\newcommand{\NOb}{N''}
\newcommand{\NOm}{\Z^{N+1}/\ku{1}}
\newcommand{\NQOm}{\Q^{N+1}/\ku{1}}
\newcommand{\pNNb}{\til{\pi}''}    %
\newcommand{\pNONb}{\pi''}    %
\newcommand{\SOb}{\Sigma''}

\newcommand{\Na}{N^{\prime\prime}} %

\newcommand{\kbc}{\A} %
\DeclareMathOperator{\Gr}{Grass} %
\DeclareMathOperator{\Sl}{SL} %
\DeclareMathOperator{\Gl}{GL} %
\def\fund{\mathfrak f} %
\def\Sym{\mathfrak S} %
\DeclareMathOperator{\id}{id}  %

\DeclareMathOperator{\CDiv}{CaDiv}

\DeclareMathOperator{\Cone}{Cone}

\newcommand{\A}{\mathbb A}
\newcommand{\C}{\mathbb C}

\newcommand{\K}{\mathbb K}
\newcommand{\N}{\mathbb N}
\newcommand{\PP}{\mathbb P}
\newcommand{\Q}{\mathbb Q}
\newcommand{\R}{\mathbb R}

\newcommand{\Z}{\mathbb Z}

\newcommand{\CO}{{\mathcal O}}

\newcommand{\CS}{{\mathcal S}}

\DeclareMathOperator{\Spec}{Spec}

\newcommand{\conv}{\operatorname{conv}} %
\DeclareMathOperator{\innt}{int}

\renewcommand{\div}{{\rm div}}
\DeclareMathOperator{\supp}{supp}

\DeclareMathOperator{\gHom}{Hom}

\newcommand{\kitem}{\begin{itemize}\vspace{-2ex}}
\newcommand{\kenditem}{\vspace{-1ex}\end{itemize}}
\newcommand{\til}[1]{{\widetilde{#1}}}
\newcommand{\surj}{\rightarrow\hspace{-0.8em}\rightarrow}

\newcommand{\ko}{\overline}

\newcommand{\ku}{\underline}

\newcommand{\kss}{\scriptscriptstyle}

\newcommand{\kbb}{{\kss \bullet}}

\DeclareMathOperator{\spann}{span}

\newtheorem{theorem}{Theorem}[section]

\newtheorem{lemma}[theorem]{Lemma}

\newtheorem{proposition}[theorem]{Proposition}

\theoremstyle{definition}
\newtheorem{definition}[theorem]{Definition}

\newtheorem{example}[theorem]{Example}

\newtheorem{recipe}[theorem]{Recipe}
\newtheorem{problem}[theorem]{Problem}

\theoremstyle{remark}
\newtheorem*{remark}{Remark}

\numberwithin{equation}{section}

\newcommand{\kst}{\,|\;}

\reversemarginpar
\newcommand{\kkk}[1]{}

\newcommand{\ktrash}[1]{{}}
\newcommand{\kpriv}[1]{{}}
\newcommand{\klabel}[1]{\label{#1}\kkk{#1}}
\definecolor{changes}{rgb}{0.1,0.65,0.03}

\newcounter{Abschnitt}[section]

\begin{document}
\title[A fansy divisor on $\ko{M}_{0,n}$]{
A fansy divisor on $\ko{M}_{0,n}$}
\author[K.~Altmann]{Klaus Altmann%
}
\address{Fachbereich Mathematik und Informatik, 
Freie Universit\"at Berlin,
Arnimalle 3, 
14195 Berlin, 
Germany}
\email{altmann@math.fu-berlin.de}
\author[G.~Hein]{Georg Hein} 
\address{Universtit\"at Duisburg-Essen,
Fachbereich 6, Mathematik, 45117 Essen, Germany}
\email{ghein@math.fu-berlin.de}

\date{}
\maketitle
\begin{abstract}
We study the relation between projective $T$-varieties and their affine cones
in the language of the so-called divisorial fans and polyhedral divisors.
As an application, we present the Grassmannian $\Gr(2,n)$
as a ``fansy divisor'' on the moduli space of stable, $n$-pointed, rational
curves $\ko{M}_{0,n}$. 
\vspace{2ex}
\end{abstract}

\section{Introduction}\label{intro}
Varieties $X$ with torus action can be described by 
divisors $\CS$ on their Chow quotients.
However, this requires the use of rather strange coefficients for $\CS$.
These coefficients have to be polyhedra if $X$ is affine, and they are polyhedral
subdivisions of some vector space in the general case. Outside of a compact region, 
these subdivisions all look alike -- they form a fan which is called the tail fan of $\CS$.
Moreover, this structure
 also gives $\CS$ its name: We will call it the ``fansy divisor'' 
of $X$.

This method was developed in
\cite{toral} and \cite{affan} -- it provides an extension of the language of toric
varieties and gives complete information about the arrangement of the $T$-orbits.
In Section \ref{fanny}, we begin with a review of the basic facts
of this theory.

We present our main result, 
how to obtain this description for projective $T$-varieties, in Theorem \ref{th-Cone}.
As an application, in Theorem \ref{th-ppfan}, we present
the Grassmannian $\Gr(2,n)$
as a fansy divisor 
$\CS=\sum_B \CS_B\otimes D_B$
on the moduli space $\ko{M}_{0,n}$ of stable rational curves
with $n$ marked points. Here, $D_B$ stands for the prime divisor consisting
of the two-component curves with a point distribution according to the 
partition $B$.
The polyhedral subdivisions
$\CS_B$ look like their common tail fan 
-- only that the origin has been replaced by 
a line segment whose direction depends on $B$.

In general, when replacing $\Gr(2,n)$ by an arbitrary generalized flag variety $G/P$,
it is difficult to obtain the Chow quotient $Y$ where the fansy divisor is supposed to
live. However, one can always determine the tail fan. It is a coarsified version
of the system of Weyl chambers associated to the semi-simple algebraic group $G$.

\section{Polyhedral and fansy divisors}\label{fanny}

Let $T$ be an affine torus over an algebraically closed field $\kk$ of 
characteristic $0$.
It gives rise to the mutually dual free abelian groups
$M:=\gHom_{\mbox{\tiny algGrp}}(T,\kk^\ast)$ and
$N:=\gHom_{\mbox{\tiny algGrp}}(\kk^\ast, T)$,
and, via $T=\Spec\kk[M]$, the torus can be recovered from them.
In \cite{toral}, we have provided a method to describe affine $T$-varieties $X$
by so-called pp-divisors on lower-dimensional algebraic varieties $Y$, 
i.e.\ by divisors having
polyhedra from $N_\Q:= N\otimes_\Z\Q$ as their coefficients.
We would like to recall some details.

\begin{definition}
\label{def-pol}
If $\sigma\subseteq N_\Q$ is a polyhedral cone, then we denote by
$\Pol(N_\Q,\sigma)$ the Grothendieck group of the semigroup 
$$
\Pol^+(N_\Q,\sigma):=\{\Delta\subseteq N_\Q\kst
\Delta = \sigma + [\mbox{compact polytope}]\}
$$
with respect to Minkowski addition. We will call 
$\tail(\Delta):=\sigma$ the tail cone of the 
elements of $\Pol(N_\Q,\sigma)$.
\vspace{-2ex}
\end{definition}

Let $Y$ be a normal and semiprojective (i.e.\ $Y\to Y_0$ is projective over
an affine $Y_0$) $\kk$-variety.
We will call $\Q$-Cartier divisors on $Y$ {\em semiample} if
a multiple of them becomes base point free.

\begin{definition}
\label{def-pp}
An element $\poldiv=\sum_i \Delta_i\otimes D_i\in 
\Pol(N_\Q,\sigma)\otimes_\Z \CaDiv(Y)$ with prime divisors $D_i$
is called a {\em pp-divisor} on $(Y,N)$ with tail cone $\sigma$
if $\Delta_i\in\Pol^+(N_\Q,\sigma)$ and if
the evaluations
$\poldiv(u):=\sum_i \min\langle \Delta_i,u\rangle D_i$
are semiample for $u\in\sigma^\vee\cap M$
and big for $u\in\innt\sigma^\vee\cap M$.
(Note that the membership $u\in\sigma^\vee:=\{u\in M_\Q\kst \langle \sigma, u\rangle \geq 0\}$
guarantees that $\min\langle \Delta_i,u\rangle > -\infty$.)
\vspace{-2ex}
\end{definition}

The common tail cone of the coefficients $\Delta_i$
will be denoted by 
$\tail(\poldiv)$.
Via $\CO_Y(\poldiv):= \oplus_{u\in\sigma^\vee\cap M} \CO_Y(\poldiv(u))$,
a pp-divisor
$\poldiv$ gives rise to the affine scheme
$X:=X(\poldiv):= \Spec \kG(Y,\CO(\poldiv))$ over $Y_0$.
The $M$-grading of its regular functions translates into an action of the torus
$T$ on $X$,
and $\tail(\poldiv)^\vee$ becomes the cone generated by the weights.
Note that the modification $\til{X}:=\til{X}(\poldiv):= \Spec_Y \CO(\poldiv)$
of $X(\poldiv)$ is a fibration over $Y$ with 
$\toric{\tail(\poldiv),N}:=\Spec\kk[\tail(\poldiv)^\vee\cap M]$ 
as general fiber.
$$ 
\xymatrix@R3ex@C5em{
\til{X} \ar[r]^{}  \ar[d]_{} &
X  \ar[d]^{} \\
Y \ar[r]^{} &
Y_0
}
$$
This construction is, in some sense, functorial.
If $\polf\in\kk(Y)\otimes_\Z N$ is a ``polyhedral, rational function'',
then $\poldiv$ and $\poldiv+\div(\polf)$ define isomorphic
graded $\CO_Y$-algebras. Up to this manipulation with
``polyhedral principal divisors'', morphisms of affine varieties
with torus action are exclusively induced by the following
morphisms $\poldiv\to\poldiv'$ of pp-divisors
(cf.\ \cite[\S8]{toral} for details):

Let $\poldiv=\sum_i\Delta_i\otimes D_i$
and $\poldiv'=\sum_j\Delta_j'\otimes D_j'$
be pp-divisors on $(Y,N)$ and $(Y',N')$ with tail cones $\sigma$ and
$\sigma'$, respectively.
If $\psi:Y \to Y'$ is such that
none of the supports of the $D_i'$ contains $\psi(Y)$,
and if $F: N\to N'$ is a linear map with
$F(\sigma)\subseteq \sigma'$,
then the relation
$$ 
\sum_i \Delta_i'\otimes \psi^\ast(D_i') =:
\psi^\ast(\poldiv')
\leq
F_\ast(\poldiv) :=
\sum_i \big(F(\Delta_i)+\sigma'\big) \otimes D_i
\vspace{-1ex}
$$
inside $\Pol(N_\Q',\sigma')\otimes_\Z\CDiv(Y)$
gives rise to an equivariant (with respect to
$T=N\otimes_\Z\kk^\ast \stackrel{F\otimes\id}{\longrightarrow}
N'\otimes_\Z\kk^\ast=T'$)
morphism
$X(\poldiv)\to X(\poldiv')$.
In particular, there are adjunction maps
$\psi^{*}(\poldiv') \to \poldiv'$
and $\poldiv \to F_{*}(\poldiv)$.
\par

In \cite{affan}, polyhedral divisors on $(Y,N)$ 
have been glued together to obtain a
description of non-affine $T$-varieties in terms of so-called divisorial fans.
To define them, we first have to broaden our previous notion of pp-divisors by
allowing the empty set $\emptyset$ as a possible coefficient for $\poldiv$:

\begin{definition}
\label{def-fanny}
Let $Y$ be a normal and semiprojective variety over $\kk$.
\vspace{-2ex}
\begin{enumerate}
\item
\label{defdef-enh}
$\poldiv=\sum_i\Delta_i\otimes D_i$ is called an
(enhanced) pp-divisor if $Z(\poldiv):=\bigcup_{\Delta_i=\emptyset} D_i$
is the support of some effective, semiample divisor on $Y$, and if
$\poldiv|_{Y\setminus Z}:=
\sum_{\Delta_i\neq\emptyset}\Delta_i\otimes D_i|_{Y\setminus Z}$
is a pp-divisor in the usual sense.
\vspace{1ex}
\item
If $\poldiv=\sum_i\Delta_i\otimes D_i$
and
$\poldiv'=\sum_i\Delta_i'\otimes D_i$ are 
(enhanced) pp-divisors on $(Y,N)$,
then we denote
$\poldiv\cap\poldiv':=\sum_i (\Delta_i\cap\Delta_i')\otimes D_i$.
It is here that empty coefficients have their natural appearance.
\vspace{1ex}
\item
\label{defdef-fanny}
A finite set $\CS=\{\poldiv^\nu\}$ is called a divisorial fan on 
$(Y,N)$ if
its elements and all their mutual intersections are 
(enhanced) pp-divisors on $Y$,
and if $X(\poldiv^\mu\cap\poldiv^\nu)\to X(\poldiv^\mu)$
is always an open embedding.
\vspace{-2ex}
\end{enumerate}
\end{definition}

Note that in (\ref{defdef-enh}), one still needs to know about
the tail cone of an enhanced pp-divisor.
The easiest possibility
to keep it as part of the data is to ask for at least one non-empty
coefficient $\Delta_i$. This is not a restriction at all, since one can always
add additional summands with the neutral element $\sigma$ as coefficient.
\vspace{0.0ex}\\
The conditions of (\ref{defdef-fanny}) can alternatively be presented in an explicit 
way. However, they then turn out to be very technical,
cf.\ \cite[5.1 and 5.3]{affan}. Hence, we will restrict ourselves 
to a special case:

First, if $\poldiv$, $\poldiv'$ are pp-divisors on $(Y,N)$ with coefficients
$\Delta_i\subseteq\Delta_i'$, then, by \cite[3.5]{affan}, 
a necessary condition
for $X(\poldiv)\to X(\poldiv')$ becoming an open embedding is that
the $\Delta_i$ are (possibly empty) faces of $\Delta_i'$,
implying that $\tail(\poldiv)\leq \tail(\poldiv')$.
In particular, the coefficients $\Delta_i^\nu$ of the elements
$\poldiv^\nu$ of a divisorial fan $\CS$ form polyhedral subdivisions
$\CS_i$ of $N_\Q$, their cells are labeled by the $\nu$'s or
$\poldiv^\nu$'s,
and their tail cones fit into the so-called tail fan
$\tail(\CS)$.

\begin{definition}
\label{def-cohfanny}
A ``{\em fansy divisor}'' is a set of pp-divisors
$\CS=\{\poldiv^\nu=\sum_i \Delta_i^v\otimes D_i\}$
on $(Y,N)$
such that for any $\mu,\nu$ there are
\vspace{-1.5ex}
\begin{enumerate}
\item
a $\umunu\in M$ and, 
numbers $c_i^{\mu\nu}$ with
$\max\langle \Delta^\mu_i, \umunu\rangle \leq c_i^{\mu\nu} \leq 
\min\langle \Delta^\nu_i, \umunu\rangle$
and
$\,\Delta^\mu_i\cap [\langle \kbb, \umunu\rangle=c_i^{\mu\nu}]=
\Delta^\nu_i\cap [\langle \kbb, \umunu\rangle=c_i^{\mu\nu}]$, and%
\vspace{1ex}
\item
an effective, semiample divisor $\Dmunu$ on
$Y\setminus Z(\poldiv^\nu)$
with 
$\,\bigcup\{D_i\kst \Delta_i^\mu\cap\Delta_i^\nu = 
\emptyset\} = \supp \Dmunu$
and $k\,\poldiv^\nu(u^{\mu\nu})- \Dmunu$ being semiample for $k\gg 0$.
\end{enumerate}
\end{definition}

\begin{proposition}[{\cite[6.9]{affan}}]
\label{prop-cohfanny}
Fansy divisors generate (via taking the finitely many mutual intersections
of the pp-divisors)
divisorial fans.
\vspace{-2ex}
\end{proposition}

\section{The tail fan of $G/P$}\label{tailGP}
Let $G$ be a semi-simple linear algebraic group; we fix
a maximal torus and a Borel subgroup $T\subseteq B \subseteq G$.
Denoting by $M$ and $N$ the mutually dual lattices of characters
and $1$-parameter subgroups of $T$, respectively, these choices
provide a system of positive roots $R^+\subseteq R$ with basis $D$.
We denote by $\Lambda_R:=\Z R:= \spann_\Z R\subseteq M$ the root lattice,
by $W:=N(T)/T=\langle s_\alpha\kst \alpha\in D\rangle$ the Weyl group,
and by $\fund:=(\N R^+)^\vee\subseteq N_\Q$ the cone of fundamental weights
of the {\em dual} root system $R^\ast$.\\
Note that $(R^+)^\vee:=(\N R^+)^\vee$ 
means the dual of the cone generated by $R^+$, but, alternatively,
we use the symbol $\alpha^\vee\in R^\ast\subseteq N$ to denote
the co-root assigned to $\alpha$.
See \cite{Springer} or \cite{Fulton}
for the basic facts concerning root systems and algebraic groups or 
Lie algebras;
see Section \ref{grassT} for the special case of $G=\Sl(n,\kk)$.

We would like to describe the $T$-action on generalized flag varieties.
Hence, we fix a subset $I\subseteq D$ and denote by $P_I\supseteq B$ 
the corresponding parabolic subgroup,
cf.\  \cite[\S 8.4]{Springer}. Its weights are 
$R^+\cup (-R_I^+)$ with
$R^+_I:= R^+\cap \Z I$ and $R_I:=R\cap\Z I$;
in particular, $P_\emptyset=B$.
Considering the subgroup
$W_I:=\langle s_\alpha\kst \alpha\in I\rangle\subseteq W$,
the set $W^I:=\{w\in W\kst w(I)>0\}$
provides nice representatives of the left
cosets of $W_I$, hence giving a bijection
$W^I\times W_I\stackrel{\sim}{\to} W$.
This splitting satisfies $\ell(w_1)+\ell(w_2)=\ell(w_1 w_2)$ with 
$\ell(w)$ referring to the length of a minimal representation
of $w\in W$ as a product of $s_\alpha$ with $\alpha\in D$.
Both $W$ and $W_I$ contain longest elements $w^0$ and $w^0_I$, respectively.
This gives rise to $w^I:=w^0 w^0_I\in W^I$.

Eventually, 
each $w\in W$ provides a set $R(w):=\{\alpha\in R^+\kst w(\alpha)<0\}$
\vspace{-0.2ex}
and \hspace{0.3em} 
$
\xymatrix@R=2.5ex{
\hspace{-0.3em}U_w:=\prod_{\alpha\in R(w)} \raisebox{-0.1ex}{$X_\alpha$}
\ar[r]_-{\sim}^-{\prod x_\alpha^{-1}}& \raisebox{0.8ex}{$\kk^{\# R(w)}$}
}
$
\vspace{-0.2ex}
where $X_\alpha$ is the 1-parameter subgroup characterized by
the isomorphisms
$x_\alpha:\kk^1\stackrel{\sim}{\to}X_{\alpha}\subseteq G$
with $x_\alpha(\alpha(t)\cdot\xi)=t\,x_\alpha(\xi)\,t^{-1}$ for $t\in T$.
Note that $\#R(w)=\ell(w)$. 
After choosing, for each $w\in W$, a representative $w\in N(T)\subseteq G$,
this leads to the Bruhat decomposition
$\,G=\bigcup_{w\in W^I} U_{w^{-1}} w P_I$.
The special choice $w:=w^I$ yields the dense, open cell,
and we denote by $\kbc:= (w^I)^{-1}U_{(w^I)^{-1}}w^I P$
a shift of it -- considered as a $T$-invariant, 
open subset $\kbc\subseteq G/P$.
Note that $\kbc$ is then an affine space.

\begin{lemma}
\klabel{lem-ueberd}
The open subsets $\{w\cdot \kbc\kst w\in W^I\}$ cover $G/P_I$. 
\vspace{-2ex}
\end{lemma}

\begin{proof}
It is sufficient to have the inclusion
$w^{-1} U_{w^{-1}}w \subseteq (w^I)^{-1}U_{(w^I)^{-1}}w^I$
for arbitrary $w\in W^I$,
and this means to check that
$v R(w^{-1})\subseteq R((w^I)^{-1})$ with $v:= w^Iw^{-1}$.%
\vspace{1.0ex}\\
If $\alpha\in R(w^{-1})$, then, the essential part is to see that
$v(\alpha)\in R^+$. If this were not true, then
$w^I$ would preserve the positivity of $-w^{-1}(\alpha)$,
i.e.\ $-w^{-1}(\alpha)\in R^+\cap (w^I)^{-1} R^+$.
However, since
$W_I(R^+\setminus R_I^+)=R^+\setminus R_I^+$ and
$w^0(R^+)=-R^+$, 
this would mean  that $-w^{-1}(\alpha)\in R^+_I$, i.e.\ that
$-\alpha\in w(R^+_I)$.
On the other hand, $w\in W^I$ means that
$w$ preserves the positivity of $I$, hence that of $R_I^+$.
Hence, the negative root $-\alpha$ cannot belong to $w(R^+_I)$.
\vspace{-2ex}
\end{proof}

Note that, since $w^I\in W^I$, the original open Bruhat cell $w^I\kbc$
belongs to the open covering provided by the previous lemma.
Moreover, the fact that $W$ acts via conjugation on $T$
implies that this canonical covering is $T$-invariant. 
Hence, it makes sense to ask for the fansy divisor $\CS$ on
some $Y$
describing the $T$-variety $G/P_I$. 
As usual, the underlying variety $Y$ 
is not uniquely determined, but there is a minimal choice --
the Chow quotient of $G/P_I$ by $T$.

\begin{proposition}
\label{prop-tail}
Assume that $I\cap (R^+\setminus R^+_I)^\bot=\emptyset$.
Then, the tail fan 
of the $T$-variety $G/P_I$ is a coarsened system of negative Weyl
chambers:
$$
\,\tail(\CS)=-W^I(W_I\fund):=\{-w\,W_I\fund\kst w\in W^I\}
$$ 
in $N_\Q$.
The lattice associated to $\CS$ is $\Lambda_R^\ast\supseteq N$;
it equals the weight lattice of the dual root system $R^\ast$.
\vspace{-2ex}
\end{proposition}

\begin{proof}
The torus $T$ acts on sets like $U_w$ via the weights $R(w)$.
Hence, the weights of the affine $T$-space $\kbc$ are
$$
(w^I)^{-1} R\big((w^I)^{-1}\big) = 
(w^I)^{-1}R^+ \cap (-R^+) =
-\big( R^+\setminus \,(w^I)^{-1} R^+\big) =
-(R^+\setminus R_{I}^+),
$$
where the last equality follows from $\,R^+\cap (w^I)^{-1} R^+ = R_I^+$
which was already used and shown in the proof of Lemma \ref{lem-ueberd}.
Thus, to get $W_I \fund$ as the tail cone of the pp-divisor describing $\kbc$,
we have to make sure that $W_I \fund=(R^+\setminus R_{I}^+)^\vee$,
which, in the special case $I=\emptyset$, 
is exactly the definition of $\fund$.
\vspace{0.5ex}\\
Since, for an element $w\in W_I$, one has 
$\langle w(\fund),\, R^+\setminus R_{I}^+ \rangle = 
\langle \fund, \,w^{-1}(R^+\setminus R_{I}^+) \rangle =
\langle \fund, \,R^+\setminus R_{I}^+ \rangle \geq 0$,
we easily obtain the inclusion $W_I\,\fund\subseteq (R^+\setminus R_{I}^+)^\vee$.
On the other hand, let $c\in N$ with $\langle c,R^+\setminus R^+_I\rangle\geq 0$.
If an $\alpha\in I\subseteq D$ fulfills $\langle c,\alpha\rangle <0$,
then, for any $\beta\in D$,
\[
\langle s_\alpha(c) , \beta\rangle = \langle c, s_\alpha(\beta)\rangle
=
\left\{\begin{array}{ll}
\langle c, -\alpha\rangle > 0 & \mbox{if $\beta=\alpha$}\\
\mbox{one out of } \langle c, D\setminus\{\alpha\}\rangle &
\mbox{if $\beta\neq\alpha$}.
\end{array}\right.
\]
Hence, $s_\alpha(c)$ has a better performance on $D$ than the original $c$
on $D$.
Since $s_\alpha\in W_I$,
induction shows that there is a $w\in W_I$ such that
$w(c)\in\fund$.
\vspace{0.5ex}\\ 
Eventually,
to describe the whole tail fan for $G/P$, we just have to apply the elements
of $W^I$ on the tail cone for $\kbc$.
However, we should also remark that the cones of $\tail(\CS)$ do not contain
any linear subspaces. This, together with the 
claim concerning the lattice structure, follows from the fact that
$\Z (R^+\setminus R_{I}^+)=\Lambda_R$.
To see this, let us take an arbitrary $\alpha\in I$. By the assumption
of the proposition, there is a
$\beta\in (R^+\setminus R_I^+)$ such that
$\langle \beta,\alpha^\vee\rangle\neq 0$.
Using \cite[Lemma 9.1.3]{Springer}, we obtain that then
at least one of the possibilities $\gamma:=\beta\pm\alpha$ belongs to $R$,
hence to $R\setminus R_I$.
Thus, we may conclude that 
$\alpha\in\spann_\Z(\beta,\gamma)\subseteq\Z (R^+\setminus R_{I}^+)$.
\end{proof}

\begin{remark}
If the assumption of Proposition \ref{prop-tail} is not satisfied, i.e.\
if $I':=I\cap (R^+\setminus R^+_I)^\bot \neq \emptyset$,
then the $T$-action admits a non-discrete kernel. This is reflected
by the fact that the weight cones are lower-dimensional or,
equivalently, that their duals $W_I\fund$ contain a common linear subspace,
namely $\spann_\Q \{\alpha^\vee\kst \alpha\in I'\}\subseteq N_\Q$.
Then, dividing this out, i.e.\ replacing $\Lambda_R^\ast$ by the corresponding
quotient, the claim of Proposition \ref{prop-tail} remains true.
\vspace{-2ex}
\end{remark}

Now, there is a standard procedure (cf.\ \cite[\S 11]{toral})
to establish the pp-divisor
$\poldiv$ for, say, the 
affine chart $\kbc\subseteq G/P_I$.

\begin{recipe}
\label{rec-pp}
Setting $\ell:=\ell(w^I)=\#(R^+\setminus R_{I}^+)$,
one has the two exact 
\vspace{-1.0ex}
sequences \hspace{0.3em}
$
\xymatrix@C=1.0em@R=0.7ex{
\hspace{-0.3em}0 \ar[r] & 
\Mb  \ar[r] &
\Z^\ell \ar[rr]^-{\deg_{\loc}} &&
\Lambda_R \ar[r] & 
0
}
$
and, more important, its dual
\[
\xymatrix@C=2.5em@R=4ex{
& N \ar@{^(->}[d]\\
0 \ar[r] & 
\Lambda_R^\ast 
\ar[rr]^-{-\langle \kbb, \,R^+\setminus R_{I}^+\rangle} 
&& \Z^\ell \ar[r]^-{\fanmap} &
\Nb \ar[r] & 
0.
}
\]
The scheme $Y'$ carrying the pp-divisor $\poldiv$ for the standard big cell
is the Chow quotient $\kbc\chQ T$. 
By \cite{chowToric}, it equals the toric variety
associated to the fan $\Sb$ arising from the coarsest common refinement
of the $\fanmap_\Q(\Q^\ell_{\geq 0}\mbox{-faces})$ in $\Nb_\Q:=\Nb\otimes_\Z\Q$.
To determine $\poldiv$ itself, we need to know, for every 
(first integral generator of a)
one-dimensional cone $c\in{\Sb}^{(1)}$, the polyhedron
$\Delta(c):=\big({\fanmap_\Q}^{\hspace{-0.2em}-1}(c)\cap\Q^\ell_{\geq 0})
- s(c)$\vspace{0.2ex},
where $s:\Nb\to\Z^\ell$ is a pre-chosen section of $\fanmap$ which does nothing
but shift all the ``fiber polytopes'' into $\Lambda_R^\ast\otimes\Q$.
Then, $\poldiv=\sum_c \Delta(c)\otimes \ko{\orb(c)}$ with
$\orb(c)\subseteq Y'$ denoting the 1-codimensional $T$-orbit corresponding to
the ray $c\in{\Sb}^{(1)}$.
\vspace{-2ex}
\end{recipe}

However, there are two problems. First, every chart $w\kbc$ leads to another
Chow quotient $Y'_w= w\kbc\chQ T$. They are all birationally equivalent,
with $w_i\kbc\chQ T$ being a blow up of $(w_i\kbc\cap w_j\kbc)\chQ T$,
and they are dominated by $Y=(G/P_I)\chQ T$. Hence,
one has to pull back all the pp-divisors $\poldiv_w$ from $Y'_w$ to $Y$
before glueing them.\\
Second, while the $Y'_w$ are at least toric with an explicitly computable,
but rather complicated fan, their common modification $Y$ is not.
Hence, one should look for those situations where $Y$ is already 
somehow known.
This leads to the case of $G/P_I=\Gr(2,n)$ where the Chow quotient has
been calculated in \cite{Kapranov}.

\section{The Grassmannian as a $T$-variety}\label{grassT}
Here, we describe the special case of the Grassmannian $\Gr(k,n)$
for $G/P_I$.
We begin with transferring notation and results from Section \ref{tailGP}
to this special case.
With $G=\Sl(n,\kk)$, the subgroups 
$T\subseteq B$ consist of the diagonal and upper triangular matrices, respectively.
In $M=\Z^n/\Z\cdot\ku{1}$, we denote by 
$L_i$ the image of the $i$-th basic vector $e_i\in\Z^n$. Then,
$D=\{\alpha_i:=L_i-L_{i+1}\kst i=1,\ldots,n-1\}$ and
$R^+=\{\alpha_{ij}:=L_i-L_j\kst i<j\}$. For each root $\alpha_{ij}$,
we have the 1-parameter family $X_{ij}=\{I_n+\xi E_{ij}\kst \xi\in\kk\}$.
Analogously, if $e^i\in\Z^n$ denotes the dual basis, then
$\alpha^\vee_{ij}:=e^i-e^j\in N=\ker(\ku{1})\subseteq\Z^n$
are the co-roots. Thus, root and weight lattice are
\[
\Lambda_R=\langle L_i-L_j\rangle \hookrightarrow
M = \langle L_i\rangle = \langle e^i-e^j \rangle ^\ast =\Lambda_W 
\hspace{1em}
(\mbox{with }
L_1+\ldots+L_n=0),
\]
satisfying $\Lambda_W/\Lambda_R\cong\Z/n\Z$.
The root system of $\Sl(n)$ is self dual. Hence,
if $\ell_i$ denotes the equivalence class of $e^i$ in $\Z^n/\Z\cdot\ku{1}$,
then
\[
\Lambda_W^\ast=\langle \ell_i-\ell_j\rangle = N \hookrightarrow 
\langle \ell_i \rangle =
\langle L_i-L_j \rangle ^\ast = \Lambda_R^\ast
\hspace{1em}
(\mbox{with }
\ell_1+\ldots+\ell_n=0)
\]
looks similar to the line above.
Since $\langle \alpha_i,  \sum_{v=1}^p \ell_v\rangle = \delta_{ip}$,
the standard Weyl chamber is 
$\fund=  \langle \ell_1, \ell_1+\ell_2,\ldots , \sum_{i=1}^{n-1}\ell_i\rangle$.
The embedding $\Sym_n\hookrightarrow\Gl(n)$ via permutation matrices
yields an isomorphism
$\Sym_n\stackrel{\sim}{\to} W$ such that the $W$-action becomes
$w(L_i)=L_{w(i)}$. 
\[
\unitlength=1.0pt
\begin{picture}(150,90)(0,-10)
\put(0.0,0.0){\line(1,0){100.00}} 
\put(0.0,0.0){\line(2,3){50.00}}
\put(100.0,0.0){\line(-2,3){50.00}} 
\put(0,0){\circle{5}}
\put(0,-5){\makebox(0,0)[ct]{$\alpha_{12}$}}
\put(50,0){\circle*{5}}
\put(50,-5){\makebox(0,0)[ct]{$\alpha_{13}$}}
\put(100,0){\circle*{5}}
\put(100,-5){\makebox(0,0)[ct]{$\alpha_{23}$}}
\put(50,25.0){\circle*{5}}
\put(50,20){\makebox(0,0)[ct]{$\alpha_{14}$}}
\put(75,37.5){\circle*{5}}
\put(80,37.5){\makebox(0,0)[lc]{$\alpha_{24}$}}
\put(50,75){\circle{5}}
\put(55,75){\makebox(0,0)[lc]{$\alpha_{34}$}}
\put(120,50){\makebox(0,0)[lt]{$R^+\setminus R_I^+\,$ for $\Gr(2,4)$}}
\end{picture}
\]
The Grassmannian $\Gr(k,n)$ is obtained via
$I:=D\setminus \{\alpha_{k}\}=
\{\alpha_1,\ldots,\widehat{\alpha_k},\ldots,\alpha_{n-1}\}$.
In particular,
$P_I=\{A\in\Sl(n)\kst A_{ij}=0 \mbox{ for } i>k,\, j\leq k\}$
and $W_I=\Sym_k\times\Sym_{n-k}\subseteq\Sym_n$. Moreover,
$W^I=\Sym_{k,n-k}:=
\{w\in\Sym_n\kst w(1)<\ldots<w(k),\,
w(k+1)<\ldots<w(n)\}$ is the set of
$(k,n-k)$-shuffles
with $\# W^I= {n \choose k}$
and
$w^I={1\,\ldots\, n \choose n\,\ldots\, 1}\cdot
{1\,\ldots\, k \choose k\,\ldots\, 1}{k+1\,\ldots\, n \choose n\,\ldots\, k+1}
=(1\,2\,\ldots\, n)^{-k}$. 
Eventually, we have
$-(R^+\setminus R_I^+) = \{L_i-L_j \kst 1\leq j\leq k < i \leq n\}$,
and the $k(n-k)$-dimensional charts $w\kbc$ turn into the usual ones with
coordinates
$a_{ij}$ where $i\notin w(1,\ldots,k)$, $j=1,\ldots,k$ and
weights $L_i-L_{w(j)}$.

\begin{proposition}
\klabel{prop-Gr-tail}
The tail fan of the $T$-variety $\Gr(k,n)$
equals 
\vspace{-0.5ex}
$$
\tail(\CS)=
\{\langle \pm\ell_1,\pm\ell_2,\ldots,\pm\ell_n\rangle\subseteq N_\Q \kst
\mbox{\rm the negative sign occurs exactly $k$ times}\}
\vspace{-0.5ex}
$$
consisting of ${n\choose k}$ cones, and the associated
lattice is $\Lambda_R^\ast=\langle\ell_i\rangle$.
\vspace{-2ex}
\end{proposition}

\begin{proof}
By Proposition \ref{prop-tail}, it remains to show that
$W_I\fund = \langle \ell_1,\ldots,\ell_k,
-\ell_{k+1},\ldots,-\ell_n\rangle$.
The left hand side is (not minimally) generated by
the elements
$\ell_J:=\sum_{j\in J}\ell_j$ with
$J\subseteq\{1,\ldots,k\}$ or $J\supseteq\{1,\ldots,k\}$.
While the first type does already fit into the claimed pattern on the right,
we use $\ell_J=-\ell_{\{1,\ldots,n\}\setminus J}$ to treat
the second.
\end{proof}

\begin{remark}
Via the projection $\Z^n\surj\Lambda_R^\ast$,
$\,e^i\mapsto \ell_i$, one obtains $\tail(\CS)$ as the image of the
${n\choose k}$ (out of $2^n$) coordinate orthants with
sign pattern $(n-k,k)$.
\vspace{-2ex}
\end{remark}

In the case of $k=2$, Kapranov has shown in \cite{Kapranov}
that $\Gr(2,n)\chQ T= \ko{M}_{0,n}$, where the latter denotes the 
moduli space of $n$-pointed, stable, rational curves.
For every partition $B=\big[B'\dcup B''=\{1,\ldots,n\}\big]$
with $\#(B'),\#(B'')\geq 2$, there is a distinguished prime divisor
$D_B$ on $\ko{M}_{0,n}$ -- it is the closure of the set
of curves with two (mutually intersecting)
components and point distribution according to $B$.

\begin{theorem}
\label{th-ppfan}
The $T$-variety $\Gr(2,n)$ corresponds to the fansy divisor
$\,\CS=\sum_B \CS_B\otimes D_B$ on 
$(\ko{M}_{0,n}, \,\Lambda_R^\ast=\langle \ell_i\rangle)$
where $\CS_B$ is the polyhedral subdivision 
arising from $\tail(\CS)=
\{\langle \pm\ell_1,\ldots,\pm\ell_n\rangle\subseteq N_\Q \kst
\mbox{\rm two negative signs}\}$
by replacing the origin
with the compact edge $C_B$ bounded by the vertices
$\frac{\#B'-1}{n-2} \,\ell_{B'}$ and
$\frac{\#B'+1-n}{n-2} \,\ell_{B'}
= \frac{\#B'-1}{n-2} \,\ell_{B'} - \ell_{B'}$.
\vspace{-2ex}
\end{theorem}

\begin{center}
\includegraphics[totalheight=4cm]{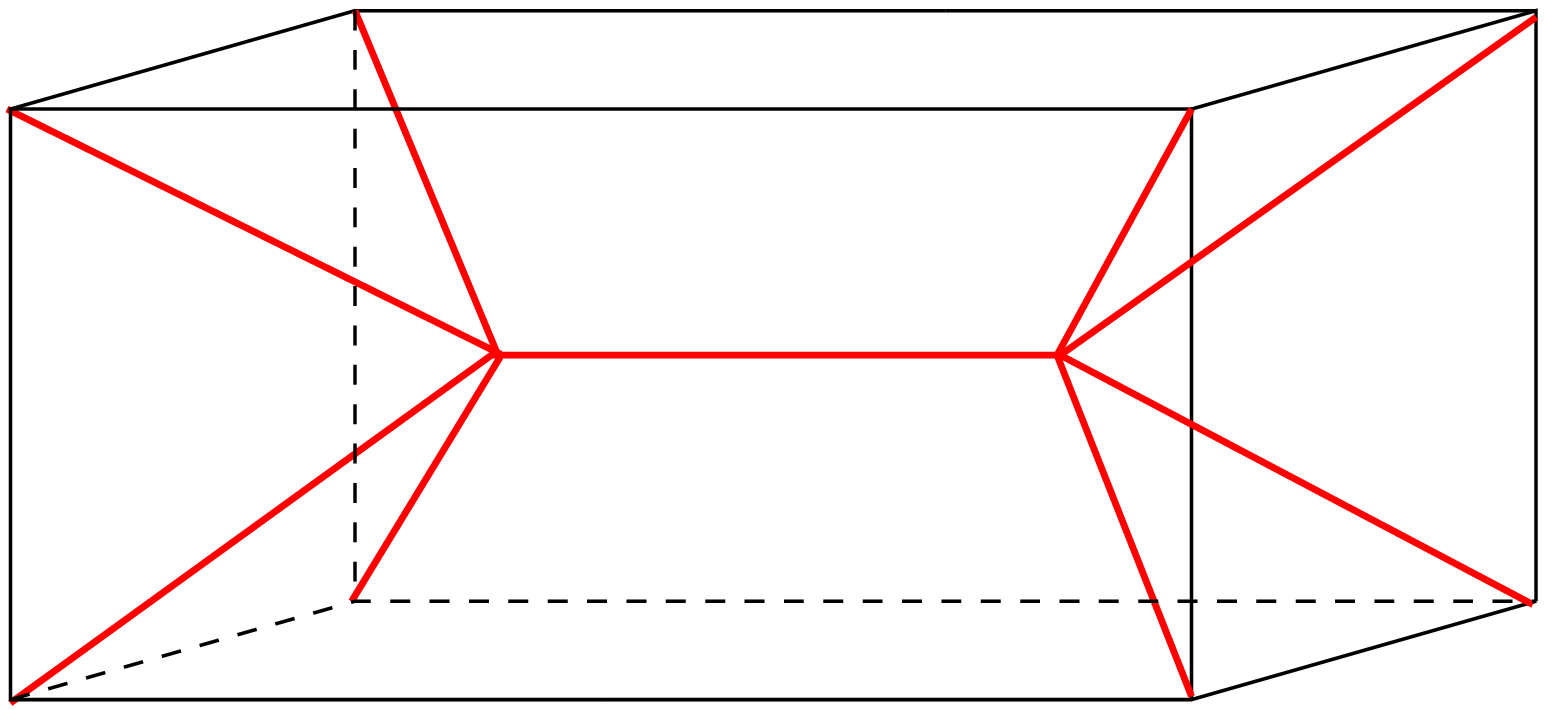}\\
{\small One of the three fansy coefficients $\CS_B
\subseteq\Q^3$ 
of $\CS$ 
on $\ko{M}_{0,4}=\PP^1$ describing $\Gr(2,4)$}
\end{center}

\begin{remark}
1) Replacing $B'$ by $B''$ does not alter the edge $C_B$. Its center is
$\frac{2\#B'-n}{2(n-2)} \,\ell_{B'}
= \frac{2\#B''-n}{2(n-2)} \,\ell_{B''}
$, and, as a vector, it equals 
$\pm\ell_{B'}=\mp\ell_{B''}$. In particular, it is the balanced
partitions $B$ that lead to edges being centered in the origin.
\vspace{0.5ex}\\
2) 
In Recipe \ref{rec-pp},
we  mentioned the
non-canonical choice of a section $s:\Nb\to\Z^\ell$.
Different choices lead to different fansy divisors --
but then their elements do only differ by a polyhedral
principal divisor as was explained in Section \ref{fanny},
right after Definition \ref{def-pp}.
Nevertheless, all ocurring polyhedral coefficients in the pp-divisors of $\CS$
are lattice polyhedra. This makes it possible to encode $\CS$ by even choosing
a {\em rational} section $s:\Nb_\Q\to\Q^\ell$. 
Among those, there is one that is canonical -- and this was used in the previous theorem.
\vspace{-2ex}
\end{remark}

The proof of Theorem \ref{th-ppfan}
does not involve a glueing of the affine charts and their corresponding
pp-divisors as mentioned at the end of Section \ref{tailGP}.
Instead, we will treat the
affine cone over $\Gr(k,n)$ with respect to its Pl\"ucker embedding. Hence,
we will take a short break and
proceed with a chapter addressing the relation between fansy divisors 
of projective
varieties and the pp-divisors of their affine cones in general.
Then, the proof of Theorem \ref{th-ppfan} will be given in 
Section \ref{affG2n}.

\section{Affine cones of projective $T$-varieties}\label{affcone}

Let $Z\subseteq\PP_\kk^N$ be a projectively normal variety 
and denote by $C(Z)\subseteq\kk^{N+1}$ its affine cone;
let them be equipped with compatible actions of 
an $n$- and an $(n+1)$-dimensional torus $\torusO$
and $\torus$, respectively. 
These actions may be described by exhibiting the degrees
of the homogeneous coordinates $z_0,\ldots,z_N$ or
of the coordinates $z_0/z_v,\ldots,z_N/z_v$ of the affine
charts $U(z_v)\subseteq\PP^N$ in the character groups
$\til{M}$ and $M$, respectively. 
Eventually, denoting by $p:\til{N}\surj N$ the projection corresponding to
$\torus\surj\torusO$, 
leads to the following commutative, mutually dual diagrams:
\[
\xymatrix{
&&0 \ar[d]& 0\ar[d]\\
0 \ar[r] & \MOb \ar@{=}[d]^-{} \ar[r] & 
\MOm \ar@{^(->}[d]
\ar[r]^-{\deg} & M \ar@{^(->}[d] \ar[r] & 0\\
0 \ar[r] & \MOb \ar[r] & \Z^{N+1} \ar[d]^-{\ku{1}} \ar[r]^-{\deg} & 
\til{M} \ar[r] \ar[d] & 0\\
&&\Z\ar@{=}[r] \ar[d]& \Z\ar[d]\\
&&0&0
}
\hspace{2em}
\xymatrix{
&0 \ar[d]& 0\ar[d]&\\
&\Z\ar@{=}[r] \ar[d]& \Z\ar[d]^-{\ku{1}}&\\
0 \ar[r] & \til{N} \ar[r]\ar@{->>}[d]^-{p} & \Z^{N+1} 
\ar@{->>}[d]^-{p}
\ar[r]_-{\pNNb} & 
\NOb \ar[r] \ar@{=}[d]^-{} \ar@/_1pc/[l]_-{s} & 0\\
0 \ar[r] & N \ar[d] \ar[r] & \NOm
\ar[d] \ar[r]^-{\pNONb} & \NOb \ar[r] & 0\\
&0&0&
}
\]
The Chow quotients do not distinguish between a projective variety 
and its affine cone -- we have
$Y:=Z\chQ \torusO = C(Z)\chQ \torus$, which is a closed subvariety of
$Y'':= \PP^{N} \chQ \torusO = \kk^{N+1}\chQ \torus$.
Similar to the local situation 
at the end of Section \ref{tailGP},
the latter equals the toric variety $\toric{\Sigma'', \NOb}$
associated to the fan $\SOb$ arising from the coarsest common refinement
of the $\pNNb_\Q(\Q^{N+1}_{\geq 0}\mbox{-faces})$ or,
equivalently, of $\pNONb_\Q(\mbox{cones of $\PP^N$-fan})$ in 
$\NOb_\Q:=\NOb\otimes_\Z\Q$.
From \cite{chowToric} and \cite{Kapranov}, we know that
$\Sigma''$ is the normal fan of the secondary polytope
$\second(\Delta)$ of
$\Delta:=\conv\{\deg z_v\kst v=0,\ldots,N\}\subseteq\til{M}_\Q$.
The faces of $\second(\Delta)$ or, equivalently,
the cones of $
\Sigma''$ correspond to the so-called regular
subdivisions of $\Delta$. This can be made explicit by
assigning to a
$\,c\in\tau\in\Sigma''$ the normal fan
\vspace{-0.5ex}
$$
\normal(\mbox{$\pNNb$}^{-1}(c)\cap \Q^{N+1}_{\geq 0})
\;\leq\; 
\normal(\ker\pNNb\cap \Q^{N+1}_{\geq 0})
\vspace{0.5ex}
$$
with ``$\leq$'' meaning ``is a subdivision of''.
Since the latter fan equals the cone over $\Delta$,
the first provides a subdivision $S(\tau)$ of the original $\Delta$,
cf.\ \cite[Lemma 2.4]{chowToric}.

\begin{lemma}
\klabel{lem-moment}
Let $\tau\in\Sigma''$ be a cone and $y\in\orb(\tau)\subseteq\toric{\Sigma'', \NOb}
=\PP^N\chQ T$.
Then, 
$y$ corresponds to a cycle $\sum_\nu \lambda_\nu Z_\nu$
with certain $T$-orbits $Z_\nu\subseteq\PP^N$, and
their images
$\mu(Z_\nu)$ under the moment map
$\mu:\PP^N\to \til{M}_\R$ yield the
subdivision $S(\tau)$.
\vspace{-2ex}
\end{lemma}

\begin{proof}
This follows directly from
\cite[Proposition 1.1]{chowToric},
its reformulation in
\cite[(0.2.10)]{Kapranov},
and the Sections
\cite[(1.2.6+7)]{Kapranov} dealing with the moment map.
\vspace{-2ex}
\end{proof}

Following Recipe \ref{rec-pp}, the pp-divisor 
$\poldiv=\sum_i\Delta_i\otimes D_i$
on $(Y,\til{N})$ that describes $C(Z)$
is built from
$\Delta_i:=\big(\mbox{$\pNNb_\Q$}^{-1}(c^i)\cap\Q^{N+1}_{\geq 0}\big)
- s(c^i)$
and $D_i:=\ko{\orb(c^i)}|_Y$
with $c^i\in\SOb^{(1)}$ browsing
through the rays of $\SOb$.
As before, $s:\NOb\to\Z^{N+1}$ denotes a section of $\pNNb$.

\begin{definition}
\label{def-boundary}
Denoting by $\{E_0,\ldots, E_N\}$ the canonical basis of
$\Z^{N+1}$ in the left diagram, we define,
for $v=0,\ldots,N$,
the faces
$\partial_v(\Delta_i):= \Delta_i\cap \big(E_v^\bot - s(c^i)\big)$
of $\Delta_i$.
\end{definition}

\begin{remark}
1) Let $\face(\Delta, E_v):=\{a\in\Delta\kst \langle a, E_v\rangle =
\min\langle \Delta, E_v \rangle\}\subseteq\til{N}_\Q$ denote the 
$\Delta$-face minimizing the
linear form $E_v$. Then,
$$
\partial_v(\Delta_i)=\left\{\begin{array}{cl}
\face(\Delta_i, E_v) & \mbox{if } \min\langle \Delta_i + s(c^i),\, E_v\rangle =0,
\mbox{ but}\\
\emptyset & \mbox{if } \min\langle \Delta_i + s(c^i),\, E_v\rangle > 0.
\end{array}\right.
$$
The union of these $\Delta_i$-faces is
$\,\bigcup_{v} \partial_v(\Delta_i) = \partial(\Delta_i)
= \big(\mbox{$\pNNb_\Q$}^{-1}(c^i)\cap\partial\Q^{N+1}_{\geq 0}\big) - s(c^i)$.
\vspace{0.5ex}\\
2) To decide whether $\partial_v(\Delta_i)$ is empty or not might be a 
cumbersome question. However, if, for $v=0,\ldots,N$, the cones
$\,\face(\tail\Delta, E_v)$ are mutually different facets of the
 cone $\tail(\Delta)$, then one {\em always} has 
$\partial_v(\Delta_i)\neq\emptyset$.
\end{remark}

\begin{example}
\label{exA3}
Let $N=2$ and consider the action of $\torus=(\kk^\ast)^2$ on $\kk^3=C(\PP^2)$
that is given by the weights 
$\deg z_0= [a,1]$, $\deg z_1= [b,1]$, $\deg z_2= [0,1]$
with relatively prime integers $a\geq b \geq 1$.
Choosing $A,B\in\Z$
with $0< A \leq b$ and $0 \leq B<a$ and
$A\,a - B\,b = 1$
yields $Y=\PP^1$, 
$\pNNb=(b,-a, a-b)$, 
$s=(-B, -A, 0)^T$, and
$E_0$, $E_1$, $E_2$ acting on $\til{N}=\Z^2$ as
$[a,1]$, $[b,1]$, and $[0,1]$, respectively.
Hence, with $c^0:=1$ and $c^{\infty}:=-1$,
\[
\poldiv \;=\;
\big(\conv\{\textstyle(\frac{B-A}{a-b},\frac{1}{a-b}),\, (\frac{A}{b},0)\}
+\kegel\big)\otimes\{0\} \;+\;
\big(\textstyle(\frac{-B}{a},0) +\kegel\big)\otimes\{\infty\}
\]
with tail cone $\kegel:=\langle(1,0), (-1,a)\rangle\subseteq\Q^2$.
We obtain the following picture of 
$p:\til{N}\surj N$ which equals the projection 
$\Z^2\surj\Z$ onto the first summand:
\vspace{1.5ex}\\
\hspace*{\fill}
\begin{picture}(160.00,150.00)(-80.00,-0.00)
\put(-100.00,110.00){\vector(0,-1){90.00}}
\put(-96.00,70.00){\makebox(0.00,0.00)[tl]{$p
$}}
\put(-30.00,90.00){\vector(-1,2){30.00}}
\put(-30.00,90.00){\line(1,-1){60.00}}
\put(30.00,30.00){\vector(1,0){40.00}}
\put(-30.00,90.00){\circle*{5.00}}
\put(-00.00,110.00){\makebox(0.00,0.00)[tc]{$(\frac{B-A}{a-b},\frac{1}{a-b})$}}
\put(30.00,30.00){\circle*{5.00}}
\put(45.00,50.00){\makebox(0.00,0.00)[tc]{$(\frac{A}{b},0)$}}
\put(0.00,0.00){\vector(1,0){80.00}}
\put(0.00,0.00){\vector(-1,0){80.00}}
\put(-55.00,5.00){\makebox(0.00,0.00)[bc]{$\Delta_0^0$}}
\put(-55.00,105.00){\makebox(0.00,0.00)[bc]{$\partial_0\Delta_0$}}
\put(-30.00,-3.00){\line(0,1){6.00}}
\put(-30.00,-5.00){\makebox(0.00,0.00)[tc]{$\frac{B-A}{a-b}$}}
\put(0.00,5.00){\makebox(0.00,0.00)[bc]{$\Delta_0^1$}}
\put(-10.00,50.00){\makebox(0.00,0.00)[bc]{$\partial_1\Delta_0$}}
\put(30.00,-3.00){\line(0,1){6.00}}
\put(30.00,-5.00){\makebox(0.00,0.00)[tc]{$\frac{A}{b}$}}
\put(50.00,15.00){\makebox(0.00,0.00)[bc]{$\partial_2\Delta_0$}}
\end{picture}
\hspace*{\fill}
\begin{picture}(110.00,150.00)(-55.00,-0.00)
\put(-10.00,30.00){\vector(-1,2){50.00}}
\put(-10.00,30.00){\vector(1,0){70.00}}
\put(-10.00,30.00){\circle*{5.00}}
\put(10.00,50.00){\makebox(0.00,0.00)[tc]{$(\frac{-B}{a},0)$}}
\put(0.00,0.00){\vector(1,0){55.00}}
\put(0.00,0.00){\vector(-1,0){55.00}}
\put(-32.00,5.00){\makebox(0.00,0.00)[bc]{$\Delta_\infty^0$}}
\put(-45.00,65.00){\makebox(0.00,0.00)[bc]{$\partial_0\Delta_\infty$}}
\put(-10.00,-3.00){\line(0,1){6.00}}
\put(-10.00,-5.00){\makebox(0.00,0.00)[tc]{$\frac{-B}{a}$}}
\put(27.00,85.00){\makebox(0.00,0.00)[bc]{$\partial_1\Delta_\infty=\emptyset$}}
\put(25.00,15.00){\makebox(0.00,0.00)[bc]{$\partial_2\Delta_\infty$}}
\end{picture}
\hspace*{\fill}
\vspace{5ex}\\
The action of the three linear forms $E_0$, $E_1$, $E_2$ is reflected by the three
different slopes of $\partial\Delta_0$.
However, it requires a closer look to realize that
$\partial_1\Delta_\infty=\emptyset$ instead of
$\,\partial_1\Delta_\infty=\{(\frac{-B}{a},0)\}$.
\end{example}

\begin{theorem}\klabel{th-Cone}
The projective $\torusO$-variety $Z\subseteq\PP_\kk^n$ is given by
the fansy divisor
$p(\partial\poldiv):=\{\sum_i \Delta_i^v\otimes D_i
\kst v=0,\ldots,N \}$
with $\Delta_i^v:= p(\partial_v\Delta_i)$.
\vspace{-2ex}
\end{theorem}

\begin{proof}
The inclusion of the two Chow quotients $Y\subseteq Y''$
may be obtained by converting $s$ into a section $\til{M}\hookrightarrow\Z^{N+1}$
and using this to interprete the original equations of 
$Z\subseteq \PP^n_\kk$ or $C(Z)\subseteq \kk^{n+1}$
as elements of $\kk[\MOb]$. This leads to a closed subscheme of
the dense torus $\Spec\kk[\MOb]\subseteq Y''$, and $Y$ arises as the normalization 
of its closure in $Y''$.
Anyway, this procedure keeps the polytopal part unchanged, and we may assume,
without loss of generality, that $Z=\PP^n_\kk$, $C(Z)=\kk^{n+1}$,
and $Y=Y''$.
\vspace{0.5ex}\\
Let $\CS=\{\poldiv^v\kst v=0,\ldots, N\}$ be the fansy divisor describing $\PP^N$
as a $T$-variety
and denote by $E^0,\ldots,E^N\in\Z^{N+1}$ the dual basis of $E_\kbb$.
Then, up to the $s$ shift, we have 
$\partial_v(\Delta_i)=\mbox{$\pNNb_\Q$}^{-1}(c^i)\cap
\langle E^0,\ldots,\widehat{E^v},\ldots,
E^N \rangle \subseteq \Q^{N+1}$.
On the other hand,
the polyhedral cone representing the $v$-th affine chart of $\PP^N_\kk$
equals $\sigma_v:=\langle p(E^0),\ldots,\widehat{p(E^v)},\ldots,
p(E^N) \rangle\subseteq\NQOm$.
Hence, up to the same $s$ shift, 
the $i$-th summand of the pp-divisor $\poldiv^v$ equals
${\pNONb_\Q}^{-1}(c^i)\cap\sigma_v \subseteq\NQOm$.
Now, the claim follows from the fact that $p:\Q^{N+1}\to \NQOm$ induces
an isomorphism from $\partial_v(\Delta_i)$ to ${\pNONb_\Q}^{-1}(c^i)\cap\sigma_v$.
\end{proof}

\begin{problem}
\label{prob-id}
If, as in Example \ref{exA3}, the assumption of the second remark after
Definition \ref{def-boundary} fails, then one can, nevertheless,
still consider the polyhedral subdivisions
$\CS_i:=p(\partial\Delta_i)$ of $N_\Q$. 
To turn the formal sum $\CS=\sum_i \CS_i\otimes D_i$ into a decent
fansy divisor $\CS$, we just need to know about the labels as explained at the end of
Section \ref{fanny}. Is there an easy, direct way to find out about them, i.e.\
to decide about the emptiness of the $\partial_v(\Delta_i)$ or their images
$\Delta_i^v$?
\end{problem}

\section{The affine cone over $\Gr(2,n)$}\label{affG2n}
This section consists of the proof of Theorem \ref{th-ppfan}.
Since we are going to use the affine cone over the Pl\"ucker embedding
$\Gr(2,n)\subseteq\PP^{{n\choose 2}-1}$,
we will specify the situation of the previous section into
\[ 
\xymatrix@C=1.3em@R=0.3ex{
0 \ar[r] & \MOb \ar[r] & \Z^{{n\choose 2}} \ar[r]^-{\deg} & 
\til{M} \ar[r] & 0\\
&& E_{ij} \ar@{|->}[r] & e_i+e_j 
}
\hspace{0.9em}
\mbox{and}
\hspace{0.9em}
\xymatrix@C=2.0em@R=0.3ex{
0 \ar[r] & \til{N} \ar[r]^-{\deg^\ast} & \Z^{{n\choose 2}} 
\ar[r]^-{\pi:=\pNNb} & \NOb \ar[r] 
& 0,\\
&& E^{ij} \ar@{|->}[r] & c^{ij} 
}
\]
where $\til{M}:=\{u\in\Z^n\kst 2|\sum_i u_i\}$.
Note that we wrote $\Z^{{n\choose 2}}$ instead of 
$\Lambda^2\Z^n$ -- the reason being that its basis elements represent the exponents
of the Pl\"ucker coordinates $z_{ij}$, i.e.\ we have
$z_{ji}=-z_{ij}$, but $E_{ji}=E_{ij}$. Moreover,
$\til{N}=\Z^n+(\frac{1}{2}+\Z)^n$, and there is a canonical surjection
$\til{N}\surj\Lambda_R^\ast$ with $e^i\mapsto\ell_i$, and
$e:=\frac{1}{2}\ku{1}$ generating the kernel.%

In Section \ref{grassT}, we considered partitions
$B=\big[B'\dcup B''=\{1,\ldots,n\}\big]$.
Extending the set of basis vectors
$E^{ij}\in\Z^{n\choose 2}$, we define for a subset
$B'\subseteq \{1,\ldots,n\}$ 
$$
E^{B'}:=\sum_{i<j,\mbox{\tiny both}\,\in B'} \hspace{-0.9em}E^{ij}
\;=\; 
{\textstyle \frac{1}{\#B'-2}}\, \sum_{b\in B'} E^{B'\setminus b},
$$
where the second equality requires that $\#B'\neq 2$.
Thus, the fact $\deg^\ast(e^i)=\sum_{j\neq i} E^{ij}$ does
upgrade to the formula
$$
\deg^\ast: e^{B'}:=\sum_{b\in B'} e^b
\;\mapsto \;
2 E^{B'} + 
\hspace{-0.8em}\sum_{i\in B'\hspace{-0.2em},\, j\in B''} \hspace{-0.8em}E^{ij}.
$$
In particular, $\deg^\ast(e^{B'}-e^{B''}) = 2(E^{B'}-E^{B''})$.
Defining $c^{B'}:=\pi(E^{B'})\in \NOb$, this implies that
$c^{B'}=
\hspace{-0.0em}\sum_{i\in B'\hspace{-0.2em},\, j\in B''}
\hspace{-0.0em}c^{ij}
=c^{B''} \hspace{-0.3em}$, 
and this element will be just called $c^B$.

\begin{lemma}
\klabel{lem-ck2fiber}
With $\,\sigma:=\til{N}_\Q\cap  \Q_{\geq 0}^{n \choose 2}=
\cone (\Delta(2,n))^\vee$, the ``positive fiber'' of $c^B$ is
\[
\pi^{-1}(c^B)\cap  \Q_{\geq 0}^{n \choose 2} 
\; = \;
\ko{E^{B'} E^{B''}}+ \sigma.
\]
In particular, $c^B$ induces the subdivision
$\Delta(2,n)=\Delta(2,n)_B'\cup\Delta(2,n)_B''$ with
the big cells
$\,\Delta(2,n)_B'=\conv\{e_i+e_j\kst i\in B'\}$ and
$\,\Delta(2,n)_B''=\conv\{e_i+e_j\kst j\in B''\}$.
\vspace{-2ex}
\end{lemma}

\begin{proof}
Let $P\in\pi^{-1}(c^B)\cap  \Q_{\geq 0}^{n \choose 2}$,
i.e.\ $P=E^{B'}+\deg^\ast(v)$
with $v=\sum_{i\in B'}\lambda_i e^i + \sum_{j\in B''} \mu_j e^j$ and
$\lambda_\kbb+\lambda_\kbb+1\geq 0$, 
$\lambda_\kbb+\mu_\kbb\geq 0$, and $\mu_\kbb+\mu_\kbb\geq 0$.
If the minimal value of $\lambda_\kbb+\lambda_\kbb$ is negative,
then we denote it by  $-c
\geq -1$.
Then, $E:= (1-c) E^{B'} + c E^{B''}$ sits on the line segment,
and $P=E+\deg^\ast(\til{v})$ with $\til{v}$ having the coordinates
$\til{\lambda_i}=\lambda_i+c/2$ and 
$\til{\mu_j}=\mu_j-c/2$,
hence satisfying $\til{\lambda_\kbb}+\til{\lambda_\kbb}\geq 0$,                      
$\til{\lambda_\kbb}+\til{\mu_\kbb}\geq 0$, and 
$\til{\mu_\kbb}+\til{\mu_\kbb}\geq 0$.
\vspace{-2ex}
\end{proof}

Since the decompositions of $\Delta(2,n)$  induced from
partitions with $\#B', \#B''\geq 2$ are proper 
and the coarsest possible, it follows that the associated elements $c^B$ are rays
in the fan $\Sigma''$,
hence, they provide divisors $\ko{\orb(c^B)}\subseteq 
\toric{\Sigma'', \NOb}=
\PP^{{n\choose 2}-1}\chQ T$.
Moreover,
these decompositions are ``matroid decompositions'' in the sense of
\cite[Definition 1.2.17]{Kapranov}. Thus, they correspond to
1-codimensional so-called Chow strata in $\Gr(2,n)\chQ T$ which are, by
Lemma \ref{lem-moment}, the restrictions of $\ko{\orb(c^B)}$ via the
Pl\"ucker embedding. 
\vspace{1ex}\\
Eventually, \cite[Theorem 4.1.8]{Kapranov} tells us that
$\Gr(2,n)\chQ T\cong \ko{M}_{0,n}$.
Under this isomorphism,
by \cite[Corollary 4.1.12]{Kapranov},
the Chow stratum $\ko{\orb(c^B)}|_{\Gr(2,n)\chQ T}$ corresponds
to the divisor $D_B\subseteq\ko{M}_{0,n}$
introduced right before Theorem \ref{th-ppfan}.

\begin{proposition}
\klabel{prop-pp-grassCone}
The affine $\torus$-variety $\cone(\Gr(2,n))$
corresponds to the pp-divisor 
$\poldiv=\sum_B\Delta_B\otimes D_B$
on $(\ko{M}_{0,n}, \Ngrass)$
with $b:=\#B'$ and
\[
\Delta_B \;=\;
\frac{b-1}{n-2} \,e^{B'} - \frac{(b-1)\,b}{2(n-2)(n-1)}\cdot \ku{1}
+\frac{1}{2} \,\ko{0\, (e^{B''}-e^{B'})}+\sigma.
\]
\vspace{-2ex}
\end{proposition}

\begin{proof}
The pp-divisor in question is the restriction of
the pp-divisor describing the affine $\torus$-space $\C^{n\choose 2}$.
Hence, for every ray $c\in\mbox{$\Sigma''$}^{(1)}$, we have to determine
the shifted ``positive fiber'' via $\pi$ -- yielding the coefficient of 
$\ko{\orb(c)}$ restricted to $\Gr(2,n)\chQ T\cong \ko{M}_{0,n}$.

First, the choice of a section $s:\NOb\to\Z^{n\choose 2}$ is equivalent to the
choice of a retraction $t:\Z^{n\choose 2}\to\Ngrass$ of $\deg^\ast$.
Moreover, according to a remark right after Theorem \ref{th-ppfan},
we prefer to preserve symmetries, hence we would rather choose a {\em rational}
retraction $t:\Q_{\geq 0}^{n \choose 2}\to\Ngrass_\Q$.
Defining
$\,t(E^{ij}):= \frac{1}{n-2} e^i + \frac{1}{n-2} e^j -
\frac{1}{(n-2)(n-1)}\, \ku{1}$,
one obtains
\[
t(E^{B'})=\frac{\#B'-1}{n-2} \,e^{B'} - \frac{(\#B'-1)\cdot\#B'}{2(n-2)(n-1)} 
\cdot\ku{1},
\]
which does nicely fit with the previous formula
$e^{B'}-e^{B''}=2(E^{B'}-E^{B''})$.
Now, we use $t$ to shift the ``positive fiber'' of $c=c^B$
from Lemma \ref{lem-ck2fiber}
into $\Ngrass_\Q$.
Since $t|_\sigma=\id$, hence $t(\sigma)=\sigma$, and
the result is
\vspace{-1ex}
\[
\renewcommand{\arraystretch}{1.5}
\Delta_B\;
\begin{array}[t]{l}
=
t\big(\pi^{-1}(c^B)\cap  \Q_{\geq 0}^{n \choose 2}\big) =
\ko{t(E^{B'}) t(E^{B''})}+ \sigma 
\\ =
t(E^{B'})+ \ko{0\, (E^{B''}-E^{B'})}+\sigma. 
\end{array}
\vspace{-2ex}
\]
On the other hand, if $c$ is not of the form $c^B$,
then Lemma \ref{lem-moment} implies that, inside
$\toric{\Sigma'', \NOb}$,
the divisor $\ko{\orb(c)}$ is disjoint to the closed subvariety $\Gr(2,n)\chQ$.
\vspace{-2ex}
\end{proof}

\definecolor{w0}{rgb}{1.0,0.0,0.0}
\definecolor{w1}{rgb}{0.0,1.0,0.0}
\definecolor{w2}{rgb}{0.0,0.0,1.0}
\hspace*{\fill}
\unitlength=1.00pt
\begin{picture}(170.0,170.0)(50.00,-15.00)
\put(0.00,0.00){\line(1,0){100.00}}
\put(100.00,0.00){\line(0,1){100.00}}
\put(100.00,100.00){\line(-1,0){100.00}}
\put(0.00,100.00){\line(0,-1){100.00}}
\multiput(0.00,0.00)(15.6,10.4){4}{\line(3,2){12.00}}
\put(100.00,0.00){\line(3,2){60.00}}
\put(100.00,100.00){\line(3,2){60.00}}
\put(0.00,100.00){\line(3,2){60.00}}
\multiput(60.00,40.00)(26,0){4}{\line(1,0){22.00}}
\put(160.00,40.00){\line(0,1){100.00}}
\put(160.00,140.00){\line(-1,0){100.00}}
\multiput(60.00,140.00)(0,-26){4}{\line(0,-1){22.00}}
\put(0.00,0.00){\circle{7.00}}
\put(-10.00,0.00){\makebox(0.00,0.00)[r]{$e^1$}}
\put(100.00,0.00){\circle{7.00}}
\put(110.00,0.00){\makebox(0.00,0.00)[l]{$e-e^2$}}
\put(0.00,100.00){\circle{7.00}}
\put(-10.00,100.00){\makebox(0.00,0.00)[r]{$e-e^4$}}
\put(100.00,100.00){\circle{7.00}}
\put(98.00,113.00){\makebox(0.00,0.00)[l]{$e^3$}}
\put(60.00,40.00){\circle{7.00}}
\put(62.50,30.00){\makebox(0.00,0.00)[c]{$e-e^3$}}
\put(160.00,40.00){\circle{7.00}}
\put(170.00,40.00){\makebox(0.00,0.00)[l]{$e^4$}}
\put(60.00,140.00){\circle{7.00}}
\put(50.00,140.00){\makebox(0.00,0.00)[r]{$e^2$}}
\put(160.00,140.00){\circle{7.00}}
\put(170.00,140.00){\makebox(0.00,0.00)[l]{$e-e^1$}}
\put(30.00,70.00){\color{w0}\line(1,0){100.00}}
\put(80.00,20.00){\color{w1}\line(0,1){100.00}}
\put(50.00,50.00){\color{w2}\line(3,2){60.00}}
\put(30.00,70.00){\color{w0}\circle*{7.00}}
\put(130.00,70.00){\color{w0}\circle*{7.00}}
\put(80.00,20.00){\color{w1}\circle*{7.00}}
\put(80.00,120.00){\color{w1}\circle*{7.00}}
\put(50.00,50.00){\color{w2}\circle*{7.00}}
\put(110.00,90.00){\color{w2}\circle*{7.00}}
\put(230.00,100.00){\makebox(0.00,0.00)[l]{$\Cone(G(2,4))$}}
\put(80.00,70.00){\circle*{4.00}}
\put(78.00,77.70){\makebox(0.00,0.00)[r]{$\frac{1}{2}e$}}
\put(170.00,25.00){\makebox(0.00,0.00)[l]{\begin{tabular}[t]{r@{\;}c@{\;}l}
Cube& = & crosscut of $\sigma$\\ &&in height 1\\
Central edges&=& compact part\\&&of $\Delta_B+\frac{1}{3}e$\end{tabular}}}
\end{picture}
\hspace*{\fill}

Now, we can finish the proof of
Theorem \ref{th-ppfan}.
According to Theorem \ref{th-Cone},
it remains to project the facets of the polyhedra $\Delta_B$ along the
map $\Ngrass_\Q\to\Lambda^\ast_R$ with
$e^i\mapsto \ell_i$ and $e^{B'}\mapsto\ell_{B'}=\sum_{i\in B'} \ell_i$.
The relation $\ell_1+\ldots+\ell_n=0$ translates into
$\ell_{B''}=-\ell_{B'}$ in $\Lambda^\ast_R$.
Hence, the vertices of the compact edge of $\Delta_B$ turn into
$\frac{b-1}{n-2} \,\ell_{B'}$ and
$\frac{n-b-1}{n-2} \,\ell_{B''} = \frac{b+1-n}{n-2} \,\ell_{B'}$. Its center is
$\frac{2b-n}{2(n-2)} \,\ell_{B'}$, and, as a vector, it equals 
$\ell_B:=\pm\ell_{B'}=\mp\ell_{B''}$. The tail cones of the $\Delta_B$-facets
yield the tail fan which was already described in
Proposition \ref{prop-Gr-tail}.
\hfill$(\Box)$

\section{Appendix: The global/local comparison}\label{app}
As additional information, a comparison of the global situation with that of the 
affine chart $\kbc=\Gr(2,n)\cap [z_{12}\neq 0]$ might provide some further insight. 
While the Pl\"ucker embedding does not go along with accompanying maps between
the lattices $\Lambda_R^\ast$ and $\til{N}$
(it is not equivariant with respect to some map 
$T\hookrightarrow (\C^\ast)^{n\choose 2}$), the retraction behaves well.
The expression of the local coordinates in terms of the quotients
$z_{ij}/z_{12}$ leads to the following commutative diagram
which connects the exact sequences of the beginning of
Section \ref{affG2n} (first line) and those
of Recipe \ref{rec-pp} at the end of Section \ref{tailGP}
(second line with $\ell=2(n-2)$).

We have used the symbols $f^i$ ($i=1,2$) and
$g^j$ ($j=3,\ldots,n$) to denote a basis for
$\Z^2$ and $\Z^{n-2}$, respectively. The latter elements
satisfy $\sum_j g^j=0$ in $\Nb$.
Moreover, for a subset
$J\subseteq\{3,\ldots,n\}$, we define
$g^J:=\sum_{j\in J} g^j$.\\
In the local chart $\kbc$, only special partitions $B$ are ``visible'' --
the elements $1$ and $2$ are not allowed to sit in a common set $B'$ or $B''$.
This is reflected in the rightmost vertical map of the
above diagram: If $B$ is ``invisible'' in $\kbc$,
then $c^B$ maps to $0$. If $B$ does separate $1$ and $2$, then
we mean with ``$B\setminus 2$'' the set among $B'$ and $B''$
originally containing $2$, but removing it then.

Applying the functor $\toric{\kbb}$ of toric varieties
to the map $(N'',\Sigma'')\to (N',\Sigma')$
yields $Y''\to Y'$. Composing this with the Pl\"ucker embedding
$Y\hookrightarrow Y''$ yields the birational modification $\ko{M}_{0,n}=Y\to Y'$
mentioned at the end of Section \ref{tailGP}.

\[
\xymatrix@R=1.5ex{
\Z \ar@{^(->}[dr]^-{\frac{1}{2}\cdot\ku{1}} 
&e^i\ar@{|->}[r] & \sum_{\kbb\neq i} E^{i\kbb}&
\\
0 \ar[r] &  
\big(\Ngrass=\Z^n+(\frac{1}{2}+\Z)^n\big) \ar[r]^-{\deg^\ast} 
\ar[dddd]^-{\mbox{$e^i\mapsto\ell_i$}}&
\raisebox{1.0ex}{$\Z^{n \choose 2}$} \ar[r]^-{\pi\,:=\,\pNNb}
\ar[dddd]^-{\parbox{7.5em}{$e^{1j}\mapsto f^2\otimes g^j\\e^{2j}\mapsto f^1\otimes g^j\\
e^{12}\mapsto-\hspace{-0.1em}\sum \hspace{-0.1em}f^\kbb\hspace{-0.3em}\otimes\! g^\kbb$}} &
\raisebox{0.0ex}{$\Na$} \ar[r] 
\ar[dddd]^-{\mbox{$c^B\mapsto g^{B\setminus 2}$}}& 
0\\
\\ \\ \\
0 \ar[r] & \big(\Lambda_R^\ast=\langle\ell_1,\ldots,\ell_n\rangle\big) 
\ar[r]^-{\deg^\ast_{\loc}} \ar[ddl]&
\Z^2 \otimes\Z^{n-2} \ar[r]^-{\pi_{\loc}:=\fanmap} &
\big(\Z^{n-2}/\ku{1}=\Nb\big)\ar[r] 
& 0
\\
& \ell_{i=1,2} \ar@{|->}[r] &
-f^i\otimes\sum_j g^j 
\\
0
& \ell_{j\geq 3} \ar@{|->}[r] &
f^1\otimes g^j + f^2\otimes g^j \\
&& f^1\otimes g^j \ar@{|->}[r] & g^j\\
&& f^2\otimes g^j \ar@{|->}[r] & -g^j.
}
\]

\end{document}